\documentclass[11pt,leqno]{article}

\usepackage{amsmath,amsfonts,amscd,amssymb,theorem}

\long\def\comment#1\endcomment{}

\comment
\endcomment

\makeatletter
\begingroup
\gdef\th@dotted{\normalfont\itshape
  \def\@begintheorem##1##2{%
        \item[\hskip\labelsep \theorem@headerfont ##1\ ##2.]}%
\def\@opargbegintheorem##1##2##3{%
   \item[\hskip\labelsep \theorem@headerfont ##1\ ##2\ (##3).]}}
\endgroup
\makeatother

\theoremstyle{dotted}

\newtheorem{theorem}{Theorem}[section]

\newtheorem{conj}[theorem]{Conjecture}

\makeatletter
\begingroup
\gdef\th@upshape{\normalfont
  \def\@begintheorem##1##2{%
        \item[\hskip\labelsep \theorem@headerfont ##1\ ##2.]}%
\def\@opargbegintheorem##1##2##3{%
   \item[\hskip\labelsep \theorem@headerfont ##1\ ##2\ (##3).]}}
\endgroup
\makeatother

\theoremstyle{upshape}

\newtheorem{defn}[theorem]{Definition}
\newtheorem{remark}[theorem]{Remark}

\makeatletter
\renewcommand{\subsection}{\@startsection{subsection}{2}{0pt}{-3ex
plus -1ex minus -0.2ex}{-2mm plus -0pt minus
-2pt}{\normalfont\bfseries}} 
\renewcommand{\subsubsection}{\@startsection{subsubsection}{3}{0pt}{-3ex
plus -1ex minus -0.2ex}{-2mm plus -0pt minus
-2pt}{\normalfont\bfseries}} 
\makeatother

\makeatletter
\@addtoreset{equation}{section}
\makeatother

\newcommand{\cntrct}                % contraction with a vector field
{\hspace{2pt}\raisebox{1pt}{\text{$\lrcorner$}}\hspace{2pt}}

\newcommand{\proof}[1][Proof.]{\smallskip\noindent{\em #1}}
\def\endproof{\hfill\ensuremath{\square}\par\medskip}

\def\eqref#1{\thetag{\ref{#1}}}

\let\latexref=\ref
\def\ref#1{{\normalfont{\latexref{#1}}}}

\newcommand{\wt}{\widetilde}

\newcommand{\idot}{{\:\raisebox{1pt}{\text{\circle*{1.5}}}}}
\newcommand{\hdot}{{\:\raisebox{3pt}{\text{\circle*{1.5}}}}}
\renewcommand{\phi}{\varphi}

\newcommand{\Hom}{\operatorname{Hom}}

\newcommand{\RHom}{\operatorname{RHom}}

\newcommand{\id}{\operatorname{\sf id}}

\newcommand{\D}{{\cal D}}

\newcommand{\Q}{\mathbb{Q}}
\newcommand{\C}{\mathbb{C}}
\newcommand{\R}{\mathbb{R}}

\newcommand{\Spec}{\operatorname{Spec}}
\newcommand{\Aut}{\operatorname{Aut}}

\newcommand{\Z}{{\mathbb Z}}

\newcommand{\calo}{\mathcal{O}}

\newcommand{\mhs}{\operatorname{\text -MHS}}
\newcommand{\hs}{\operatorname{\text -HS}}

\newcommand{\Pp}{\mathbb{P}}

\newcommand{\ch}{\operatorname{ch}}

\newcommand{\M}{\mathcal{M}}

\title{Beilinson conjecture for finite-dimensional associative algebras}

\author{D. Kaledin\thanks{Partially supported by RFBR grant
    12-01-33024, Russian Federation government grant,
    ag. 11.G34.31.0023, and the Dynasty Foundation award.}}

\begin{document}

\maketitle

\tableofcontents

\section*{Introduction}

The phrase ``Beilinson conjectures'', or, more precisely,
``Beilinson's conjectures on special values of $L$-functions''
refers to a set of conjectures formulated by A. Beilinson in
\cite{beil}. The subject is, roughly speaking, relations between
algebraic $K$-groups of a smooth proper algebraic variety $X$ over
$\Q$ and its de Rham cohomology groups equipped with the Hodge
structure. The conjectures explain and tie together nicely several
disparate known facts and present a unified picture of great beauty
and great simplicity. One wants to believe in them to such an extent
that there are some papers where the author states right away that
everything takes place in a ``Beilinson World'' -- the world where
the conjectures are true. However, the conjectures are notoriously
difficult: in the 30 years that passed since \cite{beil}, we are not
much closer to knowing whether we live in a Beilinson World or not.

Formally, the eventual goal of the conjectures is to predict, up to
a rational multiple, the values of $L$-functions of the algebraic
variety $X$ at integer points. This follows earlier conjectures by
Deligne, where the values were predicted in some cases, and expressed
in terms of the ``period matrix'' relating natural $\Q$-structures
in the de Rham cohomology and the Betti cohomology of $X$. To get
the individual values, one has to also take account of the Hodge
filtration on de Rham cohomology; in the end, for every integer, one
produces an $\R$-vector space with two independent $\Q$-structures,
and the determinant of the corresponding transition matrix should
give the $L$-values.

However, this only works in some cases; in the general case, one
ends up with a $\R$-vector space with only one obvious
$\Q$-structure, and it is not clear how to produce a number out of
that. Beilinson's insight was to add algebraic $K$-groups to the
picture. The vector space in question is the target of a certain
natural ``regulator map'' from $K$-groups; Beilinson then
conjectures that, roughly speaking, the regulator map becomes a
isomorphism after tensoring with $\R$, so that we obtain a second
$\Q$-structure, and that the determinant of the resulting period
matrix gives the $L$-values.

In the present note, we are only concerned with the first part of
the conjectures -- namely, with the prediction that the regulator
map becomes an isomorphism after tensoring with $\R$ (this is why we
use the singular ``conjecture'' in the title). We only do two
things: we observe that the conjecture can actually be stated so
that it makes sense for a non-commutative algebraic variety, and
then we observe that the conjecture thus stated is almost trivially
true for finite dimensional associative algebras (we do not really
give a complete proof, but we do give a detailed sketch). This is
boring, but not as boring as the corresponding commutative statement
-- finite-dimensional associative algebras are not as trivial as the
commutative ones. In fact, it is the basic principle of
non-commutative algebraic geomery that every variety is
derived-affine, so that were we able to extend our observation from
algebras to DG algebras, we would be close to proving the conjecture
in full generality. Needless to say, at this point we have no idea
how to do that.

The note consists of two parts. In Section~\ref{1.sec}, we review
the relevant parts of Beilinson's conjectures, with an elegant
reformulation due to K. Kato presented in Section~\ref{1a.sec}. In
Section~\ref{2.sec}, we state the non-commutative version and show
why it is true in the finite-dimensional algebra case.

\subsection*{Acknowledgements.} This paper was presented as a talk
at the annual Miami-Cinvestav-Campinas Homological Mirror Symmetry
meeting, and written up during my stay at MSRI during the
Non-commutative Geometry program in 2013. I am very grateful to the
organizers of both events. I am also grateful to the referee for
useful suggestions and corrections.

\section{Beilinson conjecture.}\label{1.sec}

\subsection{$\R$-Hodge structures.}

Recall that an {\em $\R$-mixed Hodge structure} $V$ is a triple
$\langle V_\R,V_\C,\phi \rangle$ of an $\R$-vector space $V_\R$
equipped with an increasing ``weight'' filtration $W_\idot$, a
$\C$-vector space $V_\C$ equipped with a decreasing ``Hodge''
filtration $F^\hdot$, and an isomorphism $\phi:V_\R \otimes_\R \C
\cong V_\C$ such that $W_\idot$ and $F^\hdot$ satisfy certain
non-degenracy conditions. The category $\R\mhs$ of $\R$-mixed Hodge
structures is abelian and has homological dimension $1$. 

An {\em $\R$-mixed Hodge complex} is a triple $\langle
V^\hdot_\R,V^\hdot_\C,\phi\rangle$ of a filtered complex $\langle
V^\hdot_\R,W_\idot\rangle$ of $\R$-vector spaces, a bifiltered complex
$\langle V^\hdot_\C,W_\idot,F^\hdot \rangle$ of $\C$-vector spaces, and a
filtered quasiisomorphism 
$$
\phi:\langle V_\R^\hdot,W_\idot\rangle \otimes_\R \C \cong \langle
V^\hdot_\C,W_\idot \rangle
$$
such that in every degree $i$, the homology triple $\langle
H^i(V^\hdot_\R),H^i(V^\hdot_\C),\phi \rangle$ with the induced
filtrations is an $\R$-mixed Hodge structure.

A map between mixed Hodge complexes is a {\em quasiisomorphism} if
it induces a quasiisomorphism on the associated graded quotients
with respect to filtration. One shows that inverting
quasiisomorphisms in the category of $\R$-mixed Hodge complexes, one
obtains a triangulated category equivalent to the derived category
$\D(\R\mhs)$ (this is closely related to the fact that $\R\mhs$ has
homological dimension $1$). For any $\R$-mixed Hodge complex
$\langle V^\hdot_\R,V^\hdot_\C,\phi \rangle$,
$\Hom^\hdot(\R(0),V^\hdot)$ in $\D(\R\mhs)$ can be computed by the
total complex of a bicomplex
\begin{equation}\label{beil.co.eq}
\begin{CD}
W_0V^\hdot_\R \oplus (F^0 V^\hdot_\C \cap W_0 V^\hdot_\C) @>>> W_0
V^\hdot_\C,
\end{CD}
\end{equation}
where the horizontal differential is $\phi$ on $W_0V^\hdot_\R$ and
the embedding map on the second summand.

The category $\R\mhs$ is a symmetric tensor category. The Hodge-Tate
$\R$-Hodge structure $\R(1)$ of weight $1$ is an invertible object
in this category, and for any $V^\hdot \in \D(\R\mhs)$ and any
integer $i$, one denotes $V^\hdot(i) = V^\hdot \otimes
\R(1)^{\otimes i}$.

\subsection{Absolute Hodge cohomology.}

Assume given a smooth projective algebraic variety $X$ over
$\Q$. Let $X_{an}$ be the complex manifold underlying the algebraic
variety $X \otimes_\Q \C$. We then have comparison isomorphisms
$$
H^\hdot_{DR}(X_{an}) \cong H^\hdot(X \otimes_\Q \C) \cong H^\hdot(X)
\otimes_\Q \C
$$
between de Rham cohomology groups. Moreover, by definition, all
these groups are equipped with a Hodge filtration, and the
comparison isomorphism
$$
\phi:H^\hdot(X_{an},\R) \otimes_\R \C \cong H^\hdot(X_{an},\C) \cong
H^\hdot_{DR}(X_{an})
$$
equips the Betti cohomology groups $H^\hdot(X_{an},\R)$ with a mixed
$\R$-Hodge structure (pure of weight $i$ in degree $i$).

One can also refine the comparison isomorphism to a quasiisomorphism
$$
\phi:C^\hdot(X_{an},\R) \otimes_\R \C \cong C^\hdot_{DR}(X_{an})
$$
between the Betti and the de Rham cohomology complexes, and define
an $\R$-mixed Hodge complex
\begin{equation}\label{r-ho.co.eq}
C^\hdot(X) = \langle C^\hdot(X_{an},\R), C^\hdot_{DR}(X_{an}),\phi
\rangle,
\end{equation}
where the weight filtration $W_\idot$ is the canonical filtration,
and the Hodge filtration $F^\hdot$ on the de Rham cohomology complex
$$
C^\hdot_{DR}(X_{an}) = C^\hdot(X_{an},\Omega^\hdot_X)
$$
is induced by the stupid filtration on the de Rham complex
$\Omega^\hdot_X$. The details of how one defines $\phi$ and the
cohomology complexes are completely irrelevant, since the resulting
mixed Hodge complex is anyway quasiisomorphic to the sum of its
homology objects.

\begin{defn}\label{AH.defn}
  For any integer $j$, the {\em absolute Hodge cohomology groups}\\
  $H^\hdot_{AH}(X,\R(j))$ are given by
$$
H^\hdot_{AH}(X,\R(j)) =
\Hom^\hdot_{\D(\R\mhs)}(\R(0),C^\hdot(X_{an})(j)),
$$
where $C^\hdot(X_{an})$ is the $\R$-mixed Hodge complex of
\eqref{r-ho.co.eq}.
\end{defn}

Note that by definition, absolute Hodge cohomology groups can be
computed by the bicomplex \eqref{beil.co.eq}. We also note that
since the Hodge filtration on $C^\hdot_{DR}(X_{an})$ starts with
$F^0$, we have
\begin{equation}\label{vnsh.eq}
H^i_{AH}(X,\R(j)) = 0
\end{equation}
for $i > 2j$.

\subsection{The conjecture.}

Beilinson formulates his conjectures in terms of the associated
graded quotients of $K$-groups with respect to the
$\gamma$-filtration (``motivic cohomology'' in current
parlance). For our purposes, it will be convenient to stick to
$K$-groups themselves.

To formulate the conjecture, one needs one additional
ingredient. Consider the category $\R\mhs^\iota$ of mixed $\R$-Hodge
structures $\langle V_\R,V_\C,\phi \rangle$ that are equipped with a
additional involution $\iota:V_\R \to V_\R$ and an anticomplex
involution $\iota:V_\C \to V_\C$ preserving $F^\hdot$ and $W_\idot$
and compatible with $\phi$. This is again a symmetric tensor
category of homological dimension one. Equip the Tate object
$\R(1)$ with an involution $\iota$ acting by $-1$, and let $\R(i) =
\R(1)^{\otimes i}$, $i \in \Z$. The derived category
$\D(\R\mhs^\iota)$ can then again be obtained from an obvious
$\iota$-version of the category of mixed Hodge complexes, and for
any two such complexes $V^\hdot_1$, $V^\hdot_2$, we have
$$
\RHom^\hdot_{\D(\R\mhs^\iota)}(V^\hdot_1,V^\hdot_2) \cong
\RHom^\hdot_{\D(\R\mhs)}(V^\hdot_1,V^\hdot_2)^\iota,
$$
where in the right-hand side, we take invariants with respect to the
involution $\iota$ induced by the involutions on $V^\hdot_1$ and
$V^\hdot_2$.

Assume given a smooth projective algebraic variety $X$ over
$\Q$. Note that since $X$ is defined over $\Q \subset \R$, the
complexification $X \otimes_\Q \C$ carries a natural anticomplex
involution $\iota$, and it induces an anticomplex involution of the
underlying analytic variety $X_{an}$. Therefore the mixed Hodge
complex $C^\hdot(X_{an})$ of Definition~\ref{AH.defn} becomes an
object in $\D(\R\mhs^\iota)$, and the absolute Hodge cohomology
groups acquire a natural involution $\iota$. We note that under our
convention, $H^{2i}(\Pp^n) \cong \R(-i)$ as objects in $\R\mhs^\iota$
for $0 \leq i \leq n$.

Then Beilinson shows that there exists a canonical functorial {\em
  regulator map}
\begin{equation}\label{regu}
r:K_i(X) \otimes_\Z \R \to \bigoplus_j H^{2j-i}_{AH}(X,\R(j))^\iota
\end{equation}
for any integer $i \geq 0$, and he conjectures that it is an isomorphism
for $i \geq 2$.

By virtue of \eqref{vnsh.eq}, the conjecture can be extended to
negative $i$, and in this case it is trivially true. The situation
is different for $i=0,1$, ``critical values'' in Beilinson's
terminology. This can be easily seen in the following two examples.
\begin{itemize}
\item Let $X = \Spec F$ be the spectrum of a number field. Then by
  Borel Theorem, the ranks of higher $K$-groups of $X$ are given by
$$
\dim K_i(X) \otimes_\Z \R =
\begin{cases}
0, &\quad i=2j, j \geq 1,\\
r_2, &\quad i=4j-1, j \geq 1,\\
r_1+r_2, &\quad i=4j+1, j \geq 1,
\end{cases}
$$
where $r_1$ and $r_2$ are the number of real and complex embeddings
of the field $K$. This matches exactly the right-hand side of
\eqref{regu}. However, $K_0(X)=\Z$, and $K_1(F) = F^*$, so that
$K_1(X) \otimes_\Z \R$ is infinite-dimensional.
\item Let $X$ be an elliptic curve. Then $\dim K_0(E) \otimes_\Z \R$
  is an extremely delicate invariant predicted by the
  Birch-Swinnerton-Dyer Conjecture; it is expected to be related to
  the pole order of the $L$-function of $X$ at $1$, and it depends
  very essentially on the arithmetic properties of $X$. On the other
  hand, the right-hand side of \eqref{regu} is the same for all
  elliptic curves.
\end{itemize}
To solve the obvious problem with $K_1(F)$, Beilinson introduces a
further conjecture -- he conjectures that every smooth projective
variety $X$ over $\Q$ has a regular model $X_\Z$ projective and flat
over $\Z$. He then replaces $K_\idot(X)$ in \eqref{regu} with the
image of the natural map $K_\idot(X_\Z) \to K_\idot(X)$. In the case
$X = \Spec F$, $X_\Z$ can be taken to be the spectrum of the ring of
integers $\calo_F \subset F$, and after tensoring with $\R$, the
procedure does not affect higher $K$-groups. However, it changes
$K_1$ drastically: the Dirichlet Unit Theorem says that
$$
K_1(X_\Z) \otimes_\Z \R = \calo_F^* \otimes_\Z \R
$$
has dimension $r_1+r_2-1$. This is still different from the
right-hand side of \eqref{regu}, but the difference of dimensions is
now only $1$; incidentally, it coincides with $\dim_\R K_0(X_\Z)
\otimes_Z \R$.

To take account of the remaining discrepancies, and to make things
plausible for elliptic curves, too, Beilinson includes the so-called
``height pairing'' into the picture. We will now recall this.
However, we prefer to repackage the whole story in a way that the
author learned from K. Kato's brilliant lectures at UChicago in
2010.

\section{Kato's reformulation.}\label{1a.sec}

We again start with linear algebra. Let us denote by $\R\hs$ the
category of triples $\langle V_\R,V_\C,\phi\rangle$ of an
$\R$-vector space $V_\R$, a $\C$-vector space $V_\C$ equipped witha
decreasing filtration $F^\hdot$, and an isomorphism $\phi:V_\R
\otimes_\R \C \cong V_\C$. In other words, objects in $\R\hs$ are
``$\R$-Hodge structures without the weight filtration'' -- and
without any conditions at all on the filtration $F^\hdot$. The
category $\R\hs$ is not an abelian category; however, it has a
natual structure of an exact tensor category, so that one can define
the tensor derived category $\D(\R\hs)$. As in the mixed Hodge
structure case, we can also define an ``$\R$-Hodge complex'' as a
triple $\langle V^\hdot_\R,V^\hdot_\C,\phi\rangle$ of a complex
$V^\hdot_\R$ of $\R$-vector spaces, a complex $V^\hdot_\C$ of
filtered $\C$-vector spaces, and a quasiisomorphism $\phi:V^\hdot_\R
\otimes_\R \C \cong V^\hdot_\C$; then inverting filtered
quasiisomorphisms in the category of $\R$-Hodge complexes, we obtain
the derived category $\D(\R\hs)$. We also have the obvious
$\iota$-versions $\R\hs^\iota$, $\D(\R\hs^\iota)$ of these
categories, and we have the forgetful functors
\begin{equation}\label{fgt}
\begin{aligned}
&\R\mhs \to \R\hs, \qquad\ \ \D(\R\mhs) \to \D(\R\hs),\\
&\R\mhs^\iota \to \R\hs^\iota, \qquad \D(\R\mhs^\iota) \to
\D(\R\hs^\iota).
\end{aligned}
\end{equation}
For any $\R$-Hodge complex $V^\hdot$, $\Hom^\hdot(\R(0),V^\hdot)$ in
$\D(\R\hs)$ can be computed by the total complex of a bicomplex
\begin{equation}\label{kato.co.eq}
\begin{CD}
V^\hdot_\R \oplus F^0 V^\hdot_\C @>>> V^\hdot_\C,
\end{CD}
\end{equation}
a version of \eqref{beil.co.eq} without $W_0$, and to obtain
$\Hom^\hdot$ in $\D(\R\hs^\idot)$, one has to take
$\iota$-invariants in this bicomplex.

For any smooth proper variety $X$ over $\Q$, $C^\hdot(X_{an})$
becomes an object of $\D(\R\hs^\iota)$ by applying the forgetful
functor \eqref{fgt}. The corresponding version
$$
H^\hdot_{\D}(X,\R(j)) =
\Hom^\hdot_{\D(\R\hs)}(\R(0),C^\hdot(X_{an})(j))
$$
of absolute Hodge cohomology groups of Definition~\ref{AH.defn} is
known as {\em Deligne cohomology} of $X$, and it actually predates
the absolute Hodge cohomology. As in the absolute Hodge cohomology
case, Deligne cohomology can be computed by the bicomplex
\eqref{kato.co.eq}, but in this case, one can actually take the cone
of the horizontal differential locally on $X_{an}$, that is, before
taking the cohomology complex $C^\hdot(X_{an},-)$. Then one gets actual
complexes $\R(j)$ of sheaves on $X_{an}$, and one has
$$
H^\hdot_{\D}(X,\R(j)) = H^\hdot(X_{an},\R(j)).
$$
It is in this form that the Deligne cohomology groups were originally
introduced by Deligne.

The forgetful functors \eqref{fgt} provide natural maps
$$
H^\hdot_{AH}(X,-) \to H^\hdot_{\D}(X,-),
$$
so that the regulator map \eqref{regu} induces a regulator map
\begin{equation}\label{regu.kato}
r:K_i(X_\Z) \otimes_\Z \R \to \bigoplus_j H^{2j-i}_{\D}(X,\R(j))^\iota,
\end{equation}
where $X_\Z$ is a good $\Z$-model of $X$ which we tacitly assume
given. The drawback of Deligne cohomology is that this map clearly
has no chance of being an isomorphism for all $i$ -- one checks
easily that the right-hand side does {\em not} vanish for $i <
0$. This is one of the reasons Beilinson had to refine Deligne
cohomology and introduce absolute Hodge cohomology. However, Kato
handles this discrepancy in a different way.

Namely, note that the complex $C^\hdot(X_{an})$ is equipped with a
Poincare pairing
$$
C^\hdot(X_{an}) \otimes C^\hdot(X_{an}) \to \R(d)[2d],
$$
where $d = \dim X$, and both sides are considered as objects of the
tensor category $\D(\R\hs^\iota)$. If one lets
\begin{equation}\label{wt.C}
\wt{C}^\hdot(X_{an}) = \bigoplus_jC^\hdot(X_{an})(j)[2j],
\end{equation}
then this extends to a pairing
\begin{equation}\label{pair.i.eq}
\wt{C}^\hdot(X_{an})^\iota \otimes \wt{C}^\hdot(X_{an})^\iota \to
\R(i)[2i]
\end{equation}
for any particular choice of an integer $i$ in the right-hand
side. Let us take $i=1$. Then since $\D(\R\hs^\iota)$ is a tensor
category, \eqref{pair.eq} induces a pairing
\begin{equation}\label{pair.eq}
\wt{C}^\hdot_{\D}(X_{an}) \otimes \wt{C}^\hdot_{\D}(X_{an}) \to
\RHom^\hdot(\R(0),\R(1)[2]) \to \R[1],
\end{equation}
where we denote
$$
\wt{C}^\hdot_{\D}(X_{an})^\iota =
\RHom^\hdot_{\D(\R\hs^\iota)}(\R(0),\wt{\C}^\hdot(X_{an})),
$$
considered as a complex of $\R$-vector spaces, and we use the
identification
$$
\Hom^1(\R(0),\R(1)) \cong \R
$$
easily deduced from \eqref{kato.co.eq}.

We note note that $\wt{C}^\hdot_{\D}(X_{an})$ coincides, up to a
quasiisomorphism, with the target of the regulator map
\eqref{regu.kato}. Then \eqref{pair.eq} induces a natural map
\begin{equation}\label{regu.kato.dual}
r^*:\wt{C}^\hdot_{\D}(X_{an})^\iota \to (K_\idot(X_\Z) \otimes_\Z
\R)^*[1]
\end{equation}
to the complex of $\R$-vector spaces dual to $K(X_\Z) \otimes_\Z \R$
shifted by $1$. Using this map, we can now formulate the conjecture.

\begin{conj}\label{main}
For any smooth projective variety $X$ over $\Q$ with a regular model
$X_\Z$ flat and projective over $\Z$, there exists a natural
distinguished triangle
\begin{equation}\label{tria.kato}
\begin{CD}
K_\idot(X_\Z) \otimes_\Z \R @>{r}>> \wt{C}^\hdot_{\D}(X_{an})^\iota
@>{r^*}>> (K_\idot(X_\Z) \otimes_\Z \R)^*[1] @>>>
\end{CD}
\end{equation}
of complexes of $\R$-vector spaces, where $r$ is the regulator map
\eqref{regu.kato}, and $r^*$ is the dual regulator map
\eqref{regu.kato.dual}.
\end{conj}

We note that for dimension reasons, the conjectural connecting
differential
$$
\delta:(K_\idot(X_\Z) \otimes_\Z \R)^* \to K_\idot(X_\Z) \otimes_\Z
\R
$$
in the triangle \eqref{tria.kato} can only be non-trivial in degree
$0$. It is further expected to coincide with the inverse to the
non-degenerate height pairing on the kernel of the regulator map on
$K_0(X_\Z) \otimes_\Z \R$. This explain the story for elliptic
curves: $K_0(X_\Z) \otimes_\Z \R$ can be large, but the height
pairing on the codimension-one subspace of numerically trivial cycles
is non-degenerate, so that $\delta$ kills off this subspace. On the
other hand, in the case $X = \Spec F$, $X_\Z = \Spec \calo_F$, the
map $\delta$ is actually trivial. The right-hand side of
\eqref{tria.kato} vanishes in homological degrees $\geq 2$, while in
degrees $0$ and $1$ we get two exact sequences
$$
\begin{CD}
0 @>>> \R @>>> \R^{r_1+r_2} @>>> \R^{r_1+r_2-1} @>>> 0,\\
0 @>>> \R^{r_1+r_2-1} @>>> \R^{r_1+r_2} @>>> \R @>>> 0
\end{CD}
$$
of $\R$-vector spaces, in perfect agreement with Dirichlet Unit
Theorem.

\section{Non-commutative version}\label{2.sec}

\subsection{DG algebras.}

Let us now present a version of Conjecture~\ref{main} for
non-commutative algebraic varieties. To begin with, we should
explain what does a ``non-commutative algebraic variety'' mean;
following the current practice, we take it to mean a small DG
category considered up to a Morita equivalence, as in
\cite{keller}. In fact, in the main conjecture we will even restrict
our attention to DG categories with one object, that is, to DG
algebras. However, for now, let us recall that for any small DG
category $A^\hdot$ over a commutative ring $k$, one defines its {\em
  periodic cyclic homology} $HP_\idot(A^\hdot)$. This is the
homology of a complex $CP_\idot(A^\hdot)$ of $k$-modules,
well-defined up to a quasiisomorphism. We have an additional
periodicity quasiisomorphism
\begin{equation}\label{perio}
u:CP_\idot(A^\hdot) \cong CP_\idot(A^\hdot)[2].
\end{equation}
Moreover, $CP_\idot(A^\hdot)$ can be refined to a filtered complex;
the decreasing filtration $F^\hdot$ is such that
$$
F^0CP_\idot(A^\hdot) = CP^-_\idot(A^\hdot),
$$
the negative cyclic homology comlpex, and
$$
F^iCP_\idot(A^\hdot) = u^iF^0(A^\hdot)
$$
for any integer $i$. In the filtered sense, \eqref{perio} becomes a
filtered quasiisomorphism
\begin{equation}\label{perio.filt}
u:CP_\idot(A^\hdot) \cong CP_\idot(A^\hdot)(1)[2],
\end{equation}
where $(1)$ indicates renumbering of the filtration.

Up to a filtered quasiisomorphism, the complex $CP_\idot(A^\hdot)$
is Morita-invariant -- if we have an equivalence $\D(A^\hdot) \cong
\D(B^\hdot)$ between the derived category of left modules over two
small DG categories $A^\hdot$, $B^\hdot$ given by an
$A^\hdot$-$B^\hdot$-bimodule $M$, then we have a natural
identification
$$
CP_\idot(A^\hdot) \cong CP_\idot(B^\hdot),
$$
compatible with the filtration and the periodicity isomorphism
\eqref{perio.filt}. Moreover, if the base ring $k$ is a field of
characteristic $0$, then if the derived category $\D(A^\hdot)$ is
equivalent to the derived category $\D(X)$ of quasicoherent sheaves
on a smooth proper variety $X$ over $k$, and the equivalence is
again given by a bimodule, -- loosely speaking, when $A^\hdot$ is
Morita-equivalent to $X$, -- we have a natural identification
of filtered periodic complexes
$$
CP_\idot(A^\hdot) \cong \wt{C}^\hdot_{DR}(X),
$$
where we set
$$
\wt{C}^\hdot_{DR}(X) = \bigoplus_j C^\hdot_{DR}(X)(j)[2j],
$$
as in \eqref{wt.C}. Thus periodic cyclic homology classes can be
thought of as a non-commutative generalization of de Rham cohomology
classes.

To set up an $\R$-Hodge complex structure on $CP_\idot(A^\hdot)$ in
the sense of Section~\ref{1a.sec}, we also need a version of Betti
cohomology for non-commutative varieties, and this is much more
problematic. At this point, the best candidate seems to be the {\em
  semitopological $K$-theory} of B. To\"en. We will not reproduce
its definition here (for a very brief overview with further
references, see e.g. \cite[Section 8]{icm}). Let us just say that to
any small DG category $A^\hdot$ over $\C$, one associates a certain
group-like infinite loop space $\M(A^\hdot)$, and one defines the
semitopological $K$-groups $K^{st}_\idot(A^\hdot)$ as
$$
K^{st}_\idot(A^\hdot) = \pi_\idot(\M(A^\hdot)).
$$
Since $\M(A^\hdot)$ is a grouplike infinite loop space, it defines a
spectrum, and semitopological $K$-groups are the the homotopy groups
of this spectrum. Rationally, every spectrum is an Eilenberg-Maclane
spectrum, that is, a complex; therefore we can effectively assume
that we are given a complex $K^{st}_\idot(A^\hdot)_\Q =
K^{st}_\idot(A^\hdot) \otimes \Q$ of $\Q$-vector spaces.

Semitopological $K$-theory is compatible with products of DG
categories, so that in particular, $K^{st}_\idot(\C)_\Q$ is an
algebra, and $K^{st}_\idot(A^\hdot)_\Q$ is a module over this
algebra. It is easy to prove that in fact,
$$
K^{st}_\idot(\C)_\Q \cong \Q[\beta],
$$
the algebra of polynomials in one variable $\beta$ of homological
degree $2$ (see \cite[Corollary 8.4]{icm}); the element $\beta$ is a
version of the Bott periodicity generator.

For any small DG category $A^\hdot$, one can invert the Bott
periodicity and consider the complex 
$$
K^{st}_\idot(A^\hdot)_\Q(\beta^{-1}) = K^{st}_\idot(A^\hdot)_\Q
\otimes_{\Q[\beta]} \Q[\beta,\beta^{-1}].
$$
Then a natural functorial Chern character map
\begin{equation}\label{ch.bl}
K^{st}_\idot(A^\hdot)_\Q(\beta^{-1})  \otimes_\Q \C \to CP_\idot(A^\hdot)
\end{equation}
has been constructed by A. Blanc in \cite{bl}. This maps sends
$\beta$ to the periodicity endomorphism $u$ on the right-hand side,
and it is expected to an isomorphism in good cases. More precisely,
there is the following expectation.

\begin{conj}[Blanc-To\"en]\label{conj.bt}
For any smooth and proper DG category $A^\hdot$, the map $\ch$ of
\eqref{ch.bl} is a quasiisomorphism.
\end{conj}

The precise meaning of smoothness and properness for DG categories
can be found e.g. in \cite{keller}; let us just say that both
properties hold automatically when $A^\hdot$ is Morita-equivalent to
a smooth projective algebraic variety $X$ over $\C$. In this highly
non-trivial special case, the conjecture has been recently proved by
Blanc in \cite[\S 4.1]{bl1}. The general case seems to be wide open.

\subsection{The conjecture.}

To define a non-commutative version of the middle term of the
triangle \eqref{tria.kato}, we observe that for any DG category
$A^\hdot$ over $\Z$, the semitopological $K$-theory
$K^{st}_\idot(A^\hdot \otimes \C)_\Q \otimes_\Q
\C[\beta,\beta^{-1}]$ by its very definition carries a natural
$\R$-structure (and even a $\Q$-structure). One then has the
following two choices:
\begin{enumerate}
\item one can assume Blanc-To\"en conjecture, so that the triple
$$
\langle K^{st}_\idot(A^\hdot_\idot \otimes \C)_\Q \otimes_\Q
  \R[\beta,\beta^{-1}], CP_\idot(A^\hdot \otimes \C),\ch \rangle
$$
is an $\R$-Hodge complex in the sense of Section~\ref{1a.sec}, or
\item one can define a {\em weak $\R$-Hodge complex} by dropping the
  requirement that the map $\phi:V^\hdot_\R \otimes_\R \C \to
  V^\hdot_\C$ is a quasiisomorphism, and note that the bicomplex
  \eqref{kato.co.eq} still makes perfect sense even in this more
  general situation.
\end{enumerate}
Following the path of least resistance, we will make the second
choice. For any small DG category $A^\hdot$ over $\C$, we will
denote by $CP^{\D}_\idot(A^\hdot)$ the total complex of the
bicomplex
\begin{equation}\label{nc.co.eq}
\begin{CD}
CP^-_\idot(A^\hdot) \oplus (K^{st}_\idot(A^\hdot \otimes \C)_\Q
\otimes_\Q \R[\beta,\beta^{-1}]) @>{\id - \ch}>> CP_\idot(A^\hdot),
\end{CD}
\end{equation}
where $\ch$ is the Chern character map \eqref{ch.bl}.

Now, assume given a DG algebra $A^\hdot$ over $\Z$, and recall that
it defines two small DG categories:
\begin{enumerate}
\item the category $C_\idot^{pf}(A^\hdot)$ of perfect modules over
  $A^\hdot$, and
\item the category $C_\idot^{pspf}(A^\hdot)$ of pseudoperfect
  modules over $A^\hdot$ -- that is, modules that are perfect when
  considered as modules over $\Z$.
\end{enumerate}
Here $C_\idot^{pf}(A^\hdot)$ is Morita-equivalent to $A^\hdot$. The
category $C_\idot^{pspf}(A^\hdot)$ can be very different, but at
least up to a Morita equivalence, it is well-defined. The
$\Hom$-pairing then gives a functor
$$
C_\idot^{pf}(A^\hdot) \times C_\idot^{pspf}(A^\hdot) \to
C^{pf}_\idot(\Z)
$$
into the DG category of perfect complexes of abelian groups, and
taking the periodic cyclic homology and semitopological $K$-theory,
we obtain a pairing
$$
CP_\idot(C^{pf}_\idot(A^\hdot) \otimes \C) \otimes
CP_\idot(C^{pspf}_\idot(A^\hdot) \otimes \C) \to \R(i)[2i]
$$
compatible with the weak $\R$-Hodge complex structure on both
sides. Here $i$ is an arbitrary integer, as in
\eqref{pair.i.eq}. Taking $i=1$, we obtain a pairing
\begin{equation}\label{pair.nc.eq}
CP^{\D}_\idot(C^{pf}_\idot(A^\hdot) \otimes \C) \otimes
CP^{\D}_\idot(C^{pspf}_\idot(A^\hdot) \otimes \C)[1] \to \R,
\end{equation}
a version of \eqref{pair.eq}.

At the same time, on the $K$-theory side, we have algebraic
$K$-groups $K_\idot(A^\hdot)$ for any small DG category $A^\hdot$
over $\Z$ -- for their definition, see e.g. \cite[Subsection
  5.2]{keller} -- and we have a Chern character map
$$
K_\idot(A^\hdot) \otimes \Q \to CP_\idot(A^\hdot \otimes \C).
$$
This map factor both through $CP^-_\idot(A^\hdot \otimes \Q)$ (this
is classic) and through $K^{st}_\idot(A^\hdot \otimes \C)_\Q \otimes_\Q
\Q[\beta,\beta^{-1}]$ (this has been proved in \cite{bl}). Therefore
it defines a natural regulator map
$$
r:K_\idot(A^\hdot) \otimes \R \to CP^{\D}_\idot(A^\hdot \otimes \C).
$$
Moreover, since $A^\hdot$ is defined over $\Z \subset \R \subset \C$,
the weak $\R$-Hodge complex $CP_\idot(A^\hdot \otimes \C)$ has a
natural $\iota$-version, and the regulator map actually gives a map
\begin{equation}\label{regu.nc.eq}
r:K_\idot(A^\hdot) \otimes \R \to CP^{\D}_\idot(A^\hdot \otimes \C)^\iota.
\end{equation}
Since $A^\hdot$ is Morita-equivalent to $C^{pf}_\idot(A^\hdot)$, we
may replace it with $C^{pf}_\idot(A^\hdot)$ in the right-hand side
of \eqref{regu.nc.eq}. Now define $K'_\idot(A^\hdot)_\R$ by
$$
K'_\idot(A^\hdot)_\R = \left(K_\idot(C^{pspf}_\idot(A^\hdot))
\otimes \R\right)^*,
$$
the complex of $\R$-vector spaces dual to
$K_\idot(C^{pspf}_\idot(A^\hdot)) \otimes \R$. Then by virtue of the
pairing \eqref{pair.nc.eq}, the regulator map $r$ of
\eqref{regu.nc.eq} for the DG category $C^{pspf}(A^\hdot)$ induces a
natural map
\begin{equation}\label{regu.nc.dual}
r^*:CP^{\D}_\idot(A^\hdot \otimes \C)^\iota \to
K'_\idot(A^\hdot)_\R[1],
\end{equation}
a version of \eqref{regu.kato.dual}.

\begin{conj}\label{main.nc}
For any DG algebra $A^\hdot$ over $\Z$, we have a natural
distinguished triangle
\begin{equation}\label{tria.nc}
\begin{CD}
K^\hdot(A^\hdot) \otimes \R @>{r}>>
CP^{\D}_\idot(A^\hdot \otimes \C)^\iota @>{r^*}>>
K'_\idot(A^\hdot)_\R[1] @>>>,
\end{CD}
\end{equation}
where $r$ and $r^*$ are the maps \eqref{regu.nc.eq} and
\eqref{regu.nc.dual}.
\end{conj}

\begin{remark}
For a general DG algebra $A^\hdot$, the meaning of the invariant
$K'_\idot(A^\hdot)_\R$ in Conjecture~\ref{main.nc} seems very
unclear. If $A^\hdot$ is actually smooth and proper over $\Z$, then
$C^{pf}_\idot(A^\hdot)$ is Morita-equivalent to
$C^{pspf}_\idot(A^\hdot)$, and $K'_{-i}(A^\hdot)_\R$ are the dual
$\R$-vector spaces to $K_i(A^\hdot) \otimes \R$. However, in
practice, what one cares about are DG algebras smooth and proper
over $\Q$, and it is certainly too much to expect that every such DG
algebra has a smooth and proper model over $\Z$. Note that in the
commutative case, Beilinson only asks for a model $X_\Z$ which is
regular as an abstract scheme, not smooth over $\Z$. It might very
well be that as stated, Conjecture~\ref{main.nc} admits easy
counterexamples, and one needs to impose some version of regularity;
unfortunately, at this moment we have no idea what it could possibly
be.
\end{remark}

\begin{remark}
Every map $f:A^\hdot \to B^\hdot$ between DG algebras induces a pair
of adjoint functors between the DG categories of DG modules over
them. However, only one functor of the pair always sends perfect
modules to perfect modules, and the other one always sends
pseudoperfect modules to pseudoperfect modules. Therefore we have
natural DG functors
$$
f^*:C^{pf}_\idot(A^\hdot) \to C^{pf}_\idot(B^\hdot), \qquad
f^*:C^{pspf}_\idot(B^\hdot) \to C^{pspf}_\idot(A^\hdot).
$$
Since $K$-groups are functorial with respect to DG functors, this
shows that both $K_\idot(-)$ and $K'_\idot(-)_\R$ are covariantly
functorial with respect to DG algebra maps. So are the regulator
maps $r$ and $r^*$. By adjunction, the $\Hom$-pairing is compatible
with functoriality, so that the triangles \eqref{tria.nc} of
Conjecture~\ref{main} should also be covariantly functorial in
$A^\hdot$.
\end{remark}

\subsection{Finite-dimensional algebras.}

Let us now illustrate Conjecture~\ref{main.nc} in the simple
paritucular case of finite-dimensional algebras. We thus assume
given a ring $A$ of finite rank as a $\Z$-module, and consider the
additive invariants of $A$ considered as a DG algebra over $\Z$.

We first observe that Conjecture~\ref{conj.bt} is almost trivially
true for $A^\hdot=A \otimes \C$. Here is a sketch of a possible
proof.

\proof[Sketch of a possible proof of Conjecture~\ref{conj.bt} for
  $A^\hdot = A \otimes \C$.] First of all, by the same argument as
in \cite[Lemma 8.3]{icm}, the infinite loop space $\M(A \otimes \C)$
is the group completion of the space
$$
\coprod B\Aut(P),
$$
where the sum is over the isomorphism classes of finitely generated
projective $A \otimes \C$-modules $P$, $\Aut(P)$ stands for the
group of automorphisms of such a module $P$ considered as a complex
Lie group, and $B\Aut(P)$ is its classifying space. Next, note that
since $A \otimes \C$ is a finite-dimensional algebra over $\C$, it
is a nilpotent extension of a semisimple algebra $A^{ss}$. Then
$A^{ss}$ is the sum of several matrix algebras over $\C$, and since
both $K^{st}_\idot$ and $HP_\idot$ are Morita-invariant and
compatible with sums, the statement for $A^{ss}$ follows from the
corresponding statement for $\C$ itself (that is, \cite[Corollary
  8.4]{icm}). It remains to notice that on one hand, $HP_\idot$ does
not change under nilpotent extensions of algebras over a field of
characteristic $0$ by a theorem of Goodwillie \cite{Go}, and on the
other hand, sending $P$ to $P^{ss} = P \otimes_{A \otimes \C}
A^{ss}$ gives a one-to-one correspondence between isomorphism
classes of finitely generated projective modules over $A \otimes \C$
and $A^{ss}$, and for any such projective $P$, $\Aut(P)$ is a
unipotent extension of $\Aut(P^{ss})$, so that the natural map
$B\Aut(P) \to B\Aut(P^{ss})$ is a homotopy equivalence.
\endproof

We now consider Conjecture~\ref{main.nc}, and we note that it can be
handled by the same argument. Namely, changing the $\Z$-lattice $A$
in $A \otimes \Q$ does not affect any of the terms in
\eqref{tria.nc}, and doing such a change if necessary, we may assume
that $A$ is a nilpotent extension of an algebra $A^{ss}$ of the form
$$
A^{ss} = \prod_iA_i,
$$
where $A_i$ is a maximal order in an algebra $A_i \otimes \Q$ over
$\Q$, a matrix algebra over a skew-field of finite rank over
$\Q$. For $A^{ss}$, $K_\idot(A^{ss}) \otimes \R$ and
$K'_\idot(A^{ss}) \otimes \R$ are dual vector spaces, and they can
be computed by the Borel Theorem. Thus verifying the conjecture for
$A^{ss}$ is straightforward, and we only need to show that
Conjecture~\ref{main.nc} is stable under nilpotent extensions. For
this, one combines the following two observations.
\begin{enumerate}
\item After tensoring with $\Q$, the cone of the natural map
  $K_\idot(A) \to K_\idot(A^{ss})$ is identified by the Chern
  character map with the cone of the natural map $CP^-_\idot(A
  \otimes \Q) \to CP^-_\idot(A^{ss})$ -- this is a version of the
  famous Goodwillie Theorem of \cite{Go}.
\item The other terms in the bicomplex \eqref{nc.co.eq} do not
  change inder nilpotent extensions, and by Quillen Devissage
  Theorem, neither does the right-hand side $K'_\idot(A \otimes
  \C)_\R$ of the triangle \eqref{tria.nc}.
\end{enumerate}

\bigskip

\noindent
{\sc
Steklov Math Institute, Algebraic Geometry section\\
\mbox{}\hspace{30mm}and\\
Laboratory of Algebraic Geometry, NRU HSE}

\bigskip

\noindent
{\em E-mail address\/}: {\tt kaledin@mi.ras.ru}

\end{document}